\documentclass[letterpaper, 10 pt, conference]{ieeeconf}

\IEEEoverridecommandlockouts
\usepackage{subfig}
\usepackage{float}
\usepackage{transparent}
\usepackage{amsmath}
\usepackage{amssymb}
\usepackage{graphicx}
\usepackage{color}
\usepackage{balance}

\usepackage{mathtools}

\mathtoolsset{showonlyrefs}

\def\<{\leqslant}           
\def\>{\geqslant}           

\def\[{[\![}
\def\]{]\!]}

\usepackage{transparent}
\usepackage{amsmath}
\usepackage{amssymb}
\usepackage{graphicx}
\usepackage{color}
\usepackage{fullpage}

\newcommand{\cI}{{\cal I}}
\newcommand{\cJ}{{\cal J}}
\newcommand{\cK}{{\cal K}}
\newcommand{\cL}{{\cal L}}
\newcommand{\cM}{{\cal M}}

\newcommand{\cQ}{{\cal Q}}

\newcommand{\cW}{{\cal W}}
\newcommand{\cX}{{\cal X}}

\newcommand{\gamovertwo} {{\frac{\gamma^2}{2}}}
\newcommand{\gammovertwo} \gamovertwo
\newcommand{\gammuovertwo} \gamovertwo
\newcommand{\gammuepstovertwo} \gamovertwo
\newcommand{\gammubarepstovertwo} \gamovertwo

\newcommand{\noncr} {\nonumber\\}

\newcommand{\beasnum}{\begin{eqnarray}}
\newcommand{\eeasnum}{\end{eqnarray}}
\newcommand{\beas}{\begin{eqnarray*}}
\newcommand{\eeas}{\end{eqnarray*}}

\newcommand{\be}{\begin{equation}}
\newcommand{\ee}{\end{equation}}
\newcommand{\ba}{\begin{array}}
\newcommand{\ea}{\end{array}}

\newtheorem{theorem}            {Theorem}[section]

\newtheorem{lemma}              [theorem]{Lemma}
\newtheorem{sideremark}         [theorem]{Remark}
\newtheorem{sideeg}           [theorem]{Example}
\newtheorem{sideconj}           [theorem]{Conjecture}

\def\argmin                     {\mathop{\rm argmin}}

\newcommand{\qed} {\hskip 0.2em\lower 0.7ex\hbox{\vbox{\hrule
\hbox{\vrule height 1.2ex\hskip 0.4em\vrule height 1.2ex}
\hrule}}}

\newcommand {\qform}[2]{{#1}^T {#2} {#1}}
\newcommand {\qformvec}[2]{\left(\begin{array}{cc}{#1}^T & 1\end{array}\right)  {#2} \left(\begin{array}{c}{#1} \\1\end{array}\right)}

\newcommand {\braket}[2]{\langle{#1}, {#2}\rangle}

\def \qweight {Q_\eta}

\usepackage{color}

\def \figname {fig:minplusrobustest:}
\def \eqnname {eq:minplusrobustest:}
\def \thmname {thm:minplusrobustest:}
\def \secname {sec:minplusrobustest:}

\def \lemname {sec:minplusrobustest:}

\begin{document}

\title{ Min-Plus approaches and Cluster Based Pruning for  Filtering in Nonlinear Systems}

\author{Srinivas Sridharan  \thanks{S.~Sridharan is with the Department of Mechanical and Aerospace Engineering,
University of California, San Diego, CA 92093-0411, USA. This research is supported under AFOSR grant FA9550-10-1-0233.
        { { srsridharan@eng.ucsd.edu}}}%
}

\maketitle

\begin{abstract}
The design of deterministic filters can be cast as a problem of minimizing an associated cost function for an optimal control problem.  Employing  the min-plus linearity  property of the dynamic programming operator (associated with the control problem) results in a computationally feasible approach (while avoiding linearization of the system dynamics/output). This article describes the salient features  of this approach and a specific form of pruning/projection, based on clustering, which serves to facilitate the numerical efficiency of these methods.  
\end{abstract}

\section{Introduction}\label{\secname intro}
Deterministic filtering approaches \cite{fleming1997deterministic,mortensen1968,krener2003,willems2004deterministic} have been studied as an alternative to stochastic methods of filtering. These deterministic methods are 
especially appealing in cases where the disturbance statistics are not known apriori. In fact, this approach has been successfully applied across 
various domains - quantitative finance, 
 attitude estimation \cite{mahony2008nonlinear}, etc. It is interesting
to note that in the case of linear systems with Gaussian white noise disturbances/measurement noise, both these approaches yield the same solution (the Kalman filter obtained as the solution to the associated Riccati equation):
a required feature of any optimal filtering method. The design of deterministic filters proceeds by casting the filtering problem 
as a optimal control problem on the system with dynamics that are time reversed (with respect to the original system). Solving the 
optimal control problem using  the dynamic programming method gives rise to an associated partial differential equation (the Hamilton-Jacobi-Bellman  (HJB)
equation). For systems with nonlinear dynamics/output equations most estimation schemes proceed via using a local linearized approximation model. 
Unfortunately, issues arise in the use of such methods for systems with larger nonlinearities (c.f. \cite{scholte2003}). A specific class of methods that
have emerged as a promising approach to optimal controller/filter  design for nonlinear systems -  the idempotent methods. Also termed, min/max plus 
methods these approaches exploit the fact that the dynamic programming operator in the HJB equation is a linear operator  in a space 
specifically chosen for this property \cite{maslov1987new,fleming2000max}. By  constructing a basis expansion for the value function in this space (semi-field), the solution to
the optimal control/filtering problem is rendered numerically tractable for previously more difficult classes of problems.  These ideas were
introduced for filtering in nonlinear systems in  \cite{fleming2000max} where the basis used are in the (semi-convex) dual space (obtained via the Fenchel transform). In that work,  the value function
was transformed into this space via the Fenchel transform and a fixed (albeit possible countably infinite) set of basis functions were chosen in this dual space. \textit{This is in contrast to the approach herein where we use a set of convex functions that are not fixed across different time steps}.  More recently \cite{mceneaney2008cdf} developed 
a related approach using a different representation in terms of semi-convex function basis. This was termed curse of dimensionality free methods
and have been utilized in areas  such as quantum control \cite{sjmmcdc2010}, deception games \cite{mceneaney2011idempotent}. Recent work
\cite{kallapur2012min}  introduced the application of these methods in the context of deterministic filtering for nonlinear systems. Such idempotent methods have also been applied in other areas \cite{kolokol2001idempotent}. In this article
we describe the design/implementation of the min-plus technique for deterministic filter design. Of specific novelty, are the approaches used  to help
handle the growth in the number of basis elements used to represent the value function as new measurements are made and state estimates are updated. 

The outline of the article is as follows. In Sec.~\ref{\secname secproblemformulation} we introduce the problem of interest and proceed to describe
 (in Sec.~\ref{\secname stagesofminplus}) the various stages of  these methods: (i) min-plus basis expansion, (ii) recursive update of
 basis expansion in response to new information, and  (iii) pruning of these quadratic basis terms to manage the  computational burden of the grown in
 their number. These steps are then applied to an example problem in Sec.~\ref{\secname example}. We conclude in Sec.~\ref{\secname conclusions}
with an indication of future research directions.

\section{Problem Formulation} \label{\secname secproblemformulation}
Consider a system described by
\begin{align}
x_{k+1} = A(x_k) + B w_k, \qquad y_k = C(x_k) + v_k. \label{\eqnname system}
\end{align}
In order to design a filter for this system we use the following form of the dynamics (backward dynamics equation)
\begin{align}
x_{k} = \tilde{A}(x_{k+1}) + \tilde{B} w_{k+1}, \label{\eqnname filterdyn}
\end{align}
and the output equation remains unchanged. 
Given an initial state estimate $\bar{x}_0$, the filter is obtained by minimizing the cost function 
\begin{align}
	J_T(\tilde x, w_{(\cdot)}) & \triangleq
	\frac{1}{2} \|{x}_0-\bar{x}_0\|_N^2 + \\ &\frac{1}{2}  \sum_{k=0}^{T-1}( \|w_k\|_{\qweight}^2  + \|v_{k+1}\|_{R}^2)\label{\eqnname eq:cost-functional},
\end{align}
where $N, R, \qweight$ are the weights on the different terms of interest and $v_k = y_k - C(x_k)$. 
The nonlinear optimal control problem for the system in \eqref{\eqnname filterdyn} corresponds to the following optimal cost function 
\begin{align}
	V_T(\tilde x)\triangleq\inf_{w_{(\cdot)}} J_T(\tilde{x}, w_{(\cdot)}).\label{\eqnname eq:inf}
\end{align}
Applying the dynamic programming approach to this problem leads to the following relation between the value function at consecutive time steps 
 \begin{align} \label{\eqnname vrecur1}
V_{k+1}(x) = \min_{w_0} \Big\{ V_{k}(x(k-1)| x(k) = x))  \noncr
+  \frac{1}{2} \qform{w_0}{\qweight}  + \frac{1}{2}\|y - C(x)\|^2_R \Big\},
\end{align}
where $ V_{k}(x(k-1)| x(k) = x)) $ denotes the value function at time $k-1$ given a state $x$ at time $k$. 
At   $T=0$ 
\begin{align}
V_0(x) := \frac{1}{2}  \Bigg\{ \|x - \bar{x}_0\|_{N^0}^2 + \phi^0 \Bigg\}, \label{\eqnname minplusv0}
\end{align}
which can be written in the quadratic form
\begin{align}
	\bigwedge_{i \in \cI_0}\frac{1}{2}  \qformvec{x}{Q^{v,0}_i}, \label{\eqnname v0quad}
\end{align}
where
\begin{align}
\label{\eqnname defofQ}
Q^{v,0}_i 
:= \left[\begin{array}{cc}N^0_i & {L^0_i}^T \\{L^0_i} & \bar{\phi}^0_i\end{array}\right],
\end{align}
$L_1^0= -  {\bar{x}_0}^T N^0$, $\bar\phi^0_1=\bar x_0^T N^0\bar x_0+\phi^0$ and $\cI_0=\{1\}$.
From \cite{kallapur2012min} we recall the following result
\begin{theorem}\label{\thmname filtRecursive}
Assuming that:  the value function $V_0$ has the quadratic form \eqref{\eqnname v0quad};  that  there exist $\cJ, \cL$ such that 
\begin{align} 
- \braket{y} {C(x)}_R &= \bigwedge_{j \in \cJ} \qformvec{x} {|y| Q^{c,y}_j} \label{\eqnname eq:y}\\
 \|C(x)\|^2_R &:= \bigwedge_{ l \in \cL} \qformvec{x}{Q^b_l}\label{\eqnname eq:cx};
\end{align} 
and that for all $M, \tilde{M}$ there exist  $Q^a_0$, $Q^{\tilde{a}}$ such that the following  expansions hold
\begin{align}
	&A(x)^T M A(x) := \bigwedge_{a \in \cI^0_a} \qformvec{x}{Q^a_0(M)}, \noncr
	&\tilde{M} A(x) := \bigwedge_{\tilde{a} \in {\tilde{\cI}}^0_a} \qformvec{x}{Q^{\tilde{a}}_0(\tilde{M})}\label{\eqnname aquadapprox}.
\end{align}
then,  there exist $I_1, Q^{v,1}_{k}$ such that
\begin{align}
V_1(x) = \bigwedge_{k \in \cI_1} \qformvec{x}{Q^{v,1}_k}.\label{\eqnname v1qform}
\end{align}
\qed
\end{theorem}

The optimal state estimate at any time step $k$ is given by 
\begin{align}
\hat{x}^* = \argmin_{x} V_t(x). \label{\eqnname xhatstardef}
\end{align}
It is of interest to note that $\hat{x}^*$ is in fact the argmin of one of the quadratics in the expansion of $V_k (x)$. 
Hence, given a form
\begin{align}
\frac{1}{2} [x^T\,\, 1] \left(\begin{array}{cc}q_{11} & q_{12} \\q_{21} & q_{22}\end{array}\right)  \left[\begin{array}{c}x \\1\end{array}\right],
\end{align}
the minimizing $x^*$ for this quadratic is given by 
\begin{align}
x^* = -[q_{11} + q_{11}^T]^{-1} [q_{12} + q_{21}^{T}].
\label{\eqnname eq:xstar}
\end{align}
In case the states are constrained to a specific set, then the minimization in \eqref{\eqnname xhatstardef} must be performed in the permissible set of states (which would in turn change \eqref{\eqnname eq:xstar}) .

The  above result  (Thm.~\ref{\thmname filtRecursive})  thus provides a recursion which   can  be applied repeatedly to determine/update the state estimate at every time step. In the following section, we describe the various stages involved in this approach.

\section{The stages of the min-plus approach}\label{\secname stagesofminplus}
At each time instant, the  approach introduced above involves carrying out three steps:
\begin{itemize}
\item Using the current value of the output to obtain a min-plus expansion of the output  (and associated) functions \eqref{\eqnname eq:cx}, the dynamics \eqref{\eqnname aquadapprox}.
\item Utilizing the recursion described in Thm.~\ref{\thmname filtRecursive}  to obtain the new min-plus basis expansion and the state estimate.
\item Projection/ Pruning of  the min-plus basis terms used, in order to facilitate numerical computation.
\end{itemize}

We now describe each stage in greater detail.
\subsection{Min-plus expansion of the terms in the value function}\label{\secname subsecminplusexpansion}
This step involves obtaining a quadratic approximation (in a min-plus sense) to the  terms used in \eqref{\eqnname vrecur1} 
(specifically  \eqref{\eqnname eq:cx} and \eqref{\eqnname aquadapprox}). The quadratic terms thus
obtained are used  to carry out the recursion step which yields the quadratic approximation \eqref{\eqnname v1qform}
for the value function at the next time step.  A quadratic approximation of a (continuous) function $g(x)$ over a
set $\Omega$ using a set of $L$ quadratics  is performed as follows.   Note that this window $\Omega$ is
a set centered around $\hat{x}^*$ i.e., the optimal state estimate available at that time step.  Given a set
of points $\omega_k \in \Omega$, we divide the region $\Omega$ into $L$ parts $\Omega_1, \Omega_2, \ldots \Omega_L$. For
each $\ell \in \{1, 2, \ldots,  L\}$ we design the quadratic  $h_\ell(x)$ which approximates $g(x)$ 
over $\Omega_\ell$ by solving the constrained optimization problem
\begin{align}
\min_{Q \in \cQ} \Big\{ \int_{\Omega_\ell}{\| \qformvec{x}{Q} - g(x)\|\,\, dx}  \Big\},\label{\eqnname constropt1}
\end{align}
where $\cQ$ is the constrained set for the choice of quadratics defined as
\begin{align}
\cQ := \Big\{\left[\begin{array}{cc}q_{11} & q_{12} \\{q_{12}}^T & q_{22}\end{array}\right]\bigg|q_{11} \in \mathbb{R}^{n \times n},  q_{11} > 0, q_{12}\in \mathbb{R}^{n \times 1}\Big\},
\end{align}
and subject to 
\begin{align}
\min_{x \in \Omega} \Big\{ \qformvec{x}{Q} - g(x)\Big\} \geq 0.
\end{align}

If a discrete set of points in $\Omega_\ell$ are chosen in order to evaluate a discrete form of \eqref{\eqnname constropt1}, then the problem reduces to a least square optimization (sum of square errors). In specific cases there may be a simplified/efficient  formulation to this optimization problem (this will be indicated for an example in Sec.~\ref{\secname example}).

\subsection{Recursion to produce new quadratic approximation to the value function}\label{\secname subsecrecursion}
We start with a result on combining two sets of quadratics.
\begin{lemma}\label{\lemname combinequadratics}
Given two sets of quadratics $Q_j$, $j\in\cJ$ and $Q_\ell, \ell \in \cL$; For any expression of the form
\begin{align}
\bigwedge_{j \in \cJ} Q_j + \bigwedge_{\ell \in \cL} Q_\ell,
\end{align}
there exists a set $\cM:= \cJ \times \cL$ and a corresponding set of quadratics $Q_m$ defined as
\begin{align}
Q_{m(j,\ell)} := Q_j + Q_\ell, \qquad \forall \j \in \cJ, \,\,\ell \in \cL,\,\,m(i,j) \in \cM.
\end{align}
 Further the following holds
\begin{align}
\bigwedge_{m \in \cM} Q_m = \bigwedge_{j \in \cJ} Q_j + \bigwedge_{\ell \in \cL} Q_\ell.
\end{align}
\qed
\end{lemma}

This result is essential to the recursion step. Specifically, consider the following recursion equation for the value function (for further details/derivations c.f. \cite{kallapur2012min})
\begin{align}
	&V_1(x) \nonumber \\
	&= \bigwedge_{i \in \cI_0} \Bigg\{  \frac{1}{2}  \|A(x)\|_{\cW}^2 + \Big\{L^0_i (I + B w^i_l) \nonumber \\
	&+  {w^i_c}^T B^T N^0_i (I + B w^i_l) +  {w^i_c}^T \qweight w^i_l \Big\} A(x)  \nonumber \\
	&+ L^0_i B w^i_c + \frac{1}{2} \Bigg[[{w^i_c}^T B^T N^0_i B w^i_c] + \qform{(w^i_c)}{\qweight} \noncr & \quad \quad\quad \quad\quad\quad\quad+ \bar{\phi}^0_i  \Bigg] \Bigg\}
	+ \frac{1}{2}{\|y - C(x)\|}_R^2,
\label{\eqnname v1expandedform}
\end{align}
where 
\begin{align}
w_c^i := - [\qweight +  B^T N^0_i B]^{-1}\times [B^T {{L^0_i}^T} ],\noncr
w_l^i := - [\qweight +  B^T N^0_i B]^{-1}\times  [B^T N^0_i A(x)],\noncr
\cW :=  \qform{(I + B w^i_l)}{N^0_i} + \qform{w^i_l}{\qweight}.
\label{\eqnname eq:w-opt}
\end{align}
Using \eqref{\eqnname v0quad}, \eqref{\eqnname eq:cx}, \eqref{\eqnname aquadapprox} we can write \eqref{\eqnname v1expandedform} as
\begin{align}
V_1(x) 	&= \bigwedge_{i \in \cI_0}\Bigg\{\bigwedge_{a \in \cI^0_a} \qformvec{x}{Q^a_0(M^0_i)} \nonumber \\
	&+   \bigwedge_{\tilde{a} \in {\tilde{\cI}}^0_a} \qformvec{x}{Q^{\tilde{a}}_0(\tilde{M}^0_i))} \nonumber \\
	&+ \qformvec{x}{ Q^c } \nonumber \\
	&+ \bigwedge_{m \in \cM} \qformvec{x}{Q^0_m} \Bigg\}, \label{\eqnname v1longand}
\end{align}
where
\begin{align}
Q^c &:= \left(\begin{array}{ccc}0 & 0 & 0 \\0 & 0 &0 \\0 & 0 & \phi^1_i\end{array}\right),\noncr
\phi^1_i &:= \bar{\phi}^0_i + \qform{y}{R} + 2 L^0_i B w^i_c + [{w^i_c}^T B^T N^0_i B w^i_c]\noncr.
&\qquad \qquad\qquad + \qform{(w^i_c)}{\qweight}.
	\end{align}
By applying Lem.~\ref{\lemname combinequadratics} to the set of quadratics above  we generate an index set $\cI_1$ and 
a set of quadratics $Q^{v,1}_m$ such that \eqref{\eqnname v1qform} holds.

\subsection{Projection/Pruning} \label{\secname subsecpruning}
\begin{figure}
\begin{center}
\includegraphics[width= \hsize]{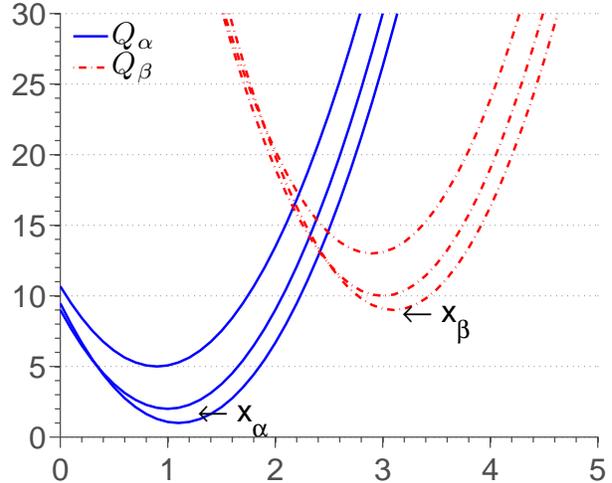}
\caption{Need for a clustering approach. Retaining only the set $Q_\alpha$ would make it harder to identify if the state shifts closer to $x_\beta$.}
\label{\figname needforclustering}
\end{center}
\end{figure}

\begin{figure}
\begin{center}
\includegraphics[width= \hsize]{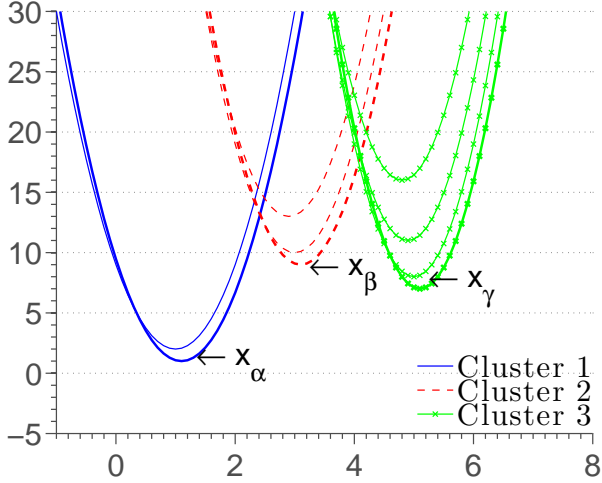}
\caption{The case when only 3 quadratics are to be retained: three clusters are formed and within each, the 
quadratic with the lowest value is selected.}
\label{\figname quadplotclustering}
\end{center}
\end{figure}

It can be seen that the construction of new quadratics leads to a growth in the number of such terms used to represent  the value function. Hence in order to facilitate numerical tractability it becomes essential to reduce the number of  such quadratics without unduly sacrificing estimation accuracy. This \lq pruning\rq\ is, in effect, a projection of an initially large set $Q$
of quadratics (min-plus basis elements)  into a lower dimensional/lower cardinality set $Q_p$ for the min-plus expansion of the value function.  
n intuitive approach to pruning in the application of the max/min plus techniques in control, involves assigning a metric which indicates 
the contribution of each quadratic to the minimization of the value function. However in the filtering case such an approach can give rise
to a non-robust filter. For instance  figure \ref{\figname needforclustering} depicts a situation where the contribution 
of the set $Q_\alpha$ of quadratics
to the minimization is much greater than that of the  
quadratics in set $Q_\beta$. However in terms of robustness, if only set $Q_\alpha$  was retained, it would cause the filter to respond more slowly to large jumps in the system state. For instance, an output corresponding to the state $\hat{x}_\beta$ would not lead to the desired state update from $\hat{x}_\alpha$ since all the quadratics around the state  $\hat{x}_\beta$ would have been pruned away.  These situations of a jump in the state,   arise in important 
applications in bimodal/multimodal oscillators, systems with binary or jump disturbances in the state.
Hence  an alternative pruning approach is required in order to increase the robustness of  the filter. One such pruning method is a clustering approach. In this technique we cluster (spatially) the  $\hat{x}^*_q$  the estimated state for each quadratic $q$ in the original set $Q$. Then,
given a requirement to retain only $P$ of these quadratics, we generate $P$ clusters and choose 
the quadratics with the best contribution metric (i.e. the ones which yield the most likely state)  
\textit{from within each cluster}. 
This proceeds as follows
\begin{enumerate}
\item  for all $q \in Q$, obtain the corresponding $\hat{x}^*_q := argmin_{x} q(x)$.  Note: due to the presence of a window of interest, 
the argmin might differ from the analytically obtained minima for the quadratic (which might lie outside the region of interest).
\item Now we cluster this set (say $\cX$) of the estimates $\hat{x}^*_q$ such that every $\hat{x} \in \cX$ belongs to a cluster $\{1, 2, \ldots, P\}$, i.e., 
the we generate  a cluster map $\cK$ which returns a cluster number for each element of $\cX$ (based on the metric chosen to separate the estimates into clusters). 
\item Within each cluster we sort and retain the quadratic that has the greatest contribution  to the minimization of the value function. Hence
the set of quadratics to be retained is defined as 
\begin{align}
Q_P := \Big \{ q \in &Q \Big| \forall \tilde{q} \in Q,\,\, \noncr &\cK(\hat{x}^*_q) = \cK(\hat{x}_{\tilde{q}}\} \Rightarrow q(\hat{x}_q) \leq \tilde{q}(\hat{x}_{\tilde{q}}) \Big \}.
\end{align}
\end{enumerate}
For a visual intuition of this pruning approach, ref Fig.~\ref{\figname quadplotclustering}.

This procedure yields the desired projection set $Q_P$ 
of quadratics, starting from the original set $Q$.

\section{Illustrative Example} \label{\secname example}

\begin{figure}
\begin{center}
\includegraphics[width= \hsize]{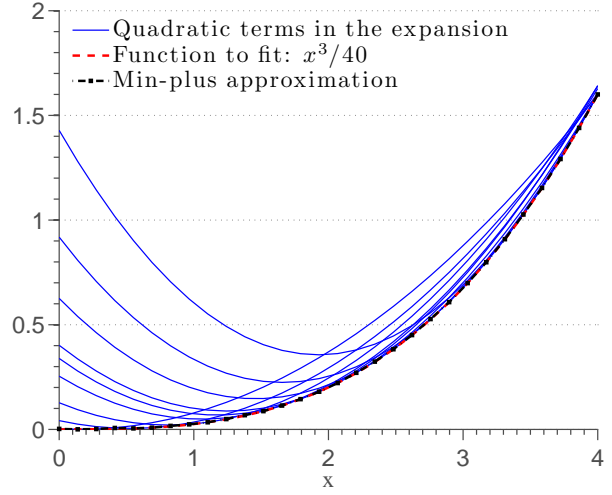}
\caption{Min-plus expansion for the output function ${x_2}^3/40$ for a window centered around $x_2 =2$.}
\label{\figname exampleminplusapprox}
\end{center}
\end{figure}

In this section we demonstrate the ideas described thus far on a two dimensional system with linear 
dynamics and a nonlinear output\footnote{The implementation code for the example in 
this paper (also useable as a template for other applications) may be
found at https://github.com/srsridharan/robustFiltering}. The continuous time representation is
\begin{align}
\frac{d}{dt} \left[\begin{array}{c}x_1 (t) \\x_2 (t) \end{array}\right]  &= \left[\begin{array}{cc}0 & 0 \\1 & 0\end{array}\right]  \left[\begin{array}{c}x_1 (t) \\x_2 (t) \end{array}\right]  + \left[\begin{array}{c}1 \\0\end{array}\right] w(t),\\
y(t) & = \frac{(x_2 (t))^3}{40}  + v(t), \label{\eqnname examplesysop}
\end{align}
where $w(\cdot)$ and $v(\cdot)$ are the process disturbance and measurement noise respectively.
Taking a sample time $\Delta t$ of $0.1$s, the discretized dynamics is
\begin{align}
\left[\begin{array}{c}x_1 (k+1) \\x_2 (k+1) \end{array}\right]  &= \left[\begin{array}{cc}1 & 0 \\0.1 & 1\end{array}\right]  \left[\begin{array}{c}x_1 (k) \\x_2 (k) \end{array}\right]  \\ & \quad + \left[\begin{array}{c}0.1 \\0\end{array}\right] \Delta w(k), \label{\eqnname discdynamicsexample}
\end{align}    
where $\Delta w(k)$ is the approximation corresponding to $w(k) \Delta t$ over the sampling time. 
As specified in Sec.~\ref{\secname stagesofminplus} the implementation of the deterministic filter proceeds along  three steps. The general  steps as modified and applied to this specific example are as follows 
to this example are as follows.\\
1. generation of the quadratic approximation: In this case, the output equation \eqref{\eqnname examplesysop} is
\begin{align}
C(x) = \frac{x^3}{40} + v(k).
\end{align}
The contrained nonlinear optimization described in Sec.~\ref{\secname subsecminplusexpansion} admits the following simplification in the current
case. Consider any window $\Omega$ with subpartitions $\Omega_k$ (recall Sec.~\ref{\secname subsecminplusexpansion}). In order to design the 
optimal quadratic approximation over such a set $\Omega_k$, we note that this problem is essentially a one dimensional regression (albeit 
with nonlinear constraints). We create a vector $x_s$ of $N$ sample points such that ${x_s}_i \in \Omega_k$ for all $i\in \{1, 2, \ldots N\}$. Now 
to fit the (continuous) output function (or its square), denoted by $g(\cdot)$, over $\Omega_k$ using a quadratic
\begin{align}
f(x) = \qformvec{x}{\left[\begin{array}{cc}a_2 & a_1/2 \\a_1/2 & a_0\end{array}\right]},
\end{align}
the optimization problem reduces to solving a least squares fitting problem to determine the optimal coefficient vector $z := [a_2, a_1, a_0]^T$ such that
\begin{align}
A z &= \left[\begin{array}{c}g({x_s}_1) \\g({x_s}_2) \\\vdots \\g({x_s}_N)\end{array}\right], \noncr
\text{where} \qquad A &:= \left[\begin{array}{ccc}({x_s}_1)^2 & ({x_s}_1) & 1 \\({x_s}_2)^2 & ({x_s}_2) & 1 \\\vdots & \vdots & \vdots \\({x_s}_N)^2 & ({x_s}_N) & 1\end{array}\right].
\end{align}
The additional constraints are that $a_2 \geq 0$ (for convexity) and that $f(x) \geq g(x), \,\,\forall x \in \Omega$. Solving this 
yields the desired quadratic bases.  The min-plus basis obained by such a quadratic approximation over a window $\Omega$ centered around $x_2 = 2$ is
as shown in Fig.~\ref{\figname exampleminplusapprox}.

2. The recursion to obtain the next set of quadratics (i.e. the min-plus basis expansion of the value function during the next time step) 
involves taking the sum of the various quadratic terms in \eqref{\eqnname v1expandedform}. 
This is generated by taking the pairwise sum of the various quadratics available.

\begin{figure}
\begin{center}
\includegraphics[width= \hsize]{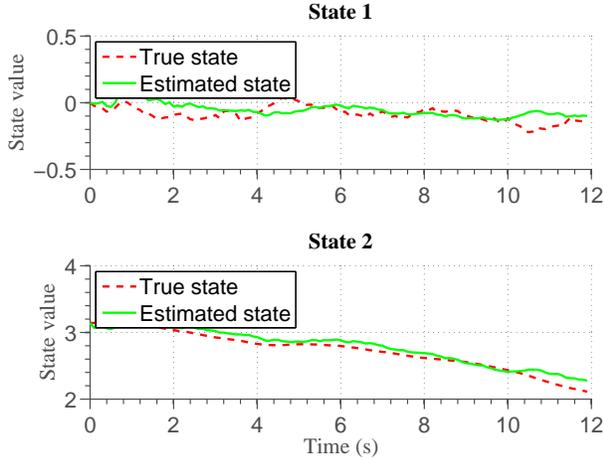}
\caption{State filtering}
\label{\figname examplesimstatefilt}
\end{center}
\end{figure}

\begin{figure}
\begin{center}
\includegraphics[width= \hsize]{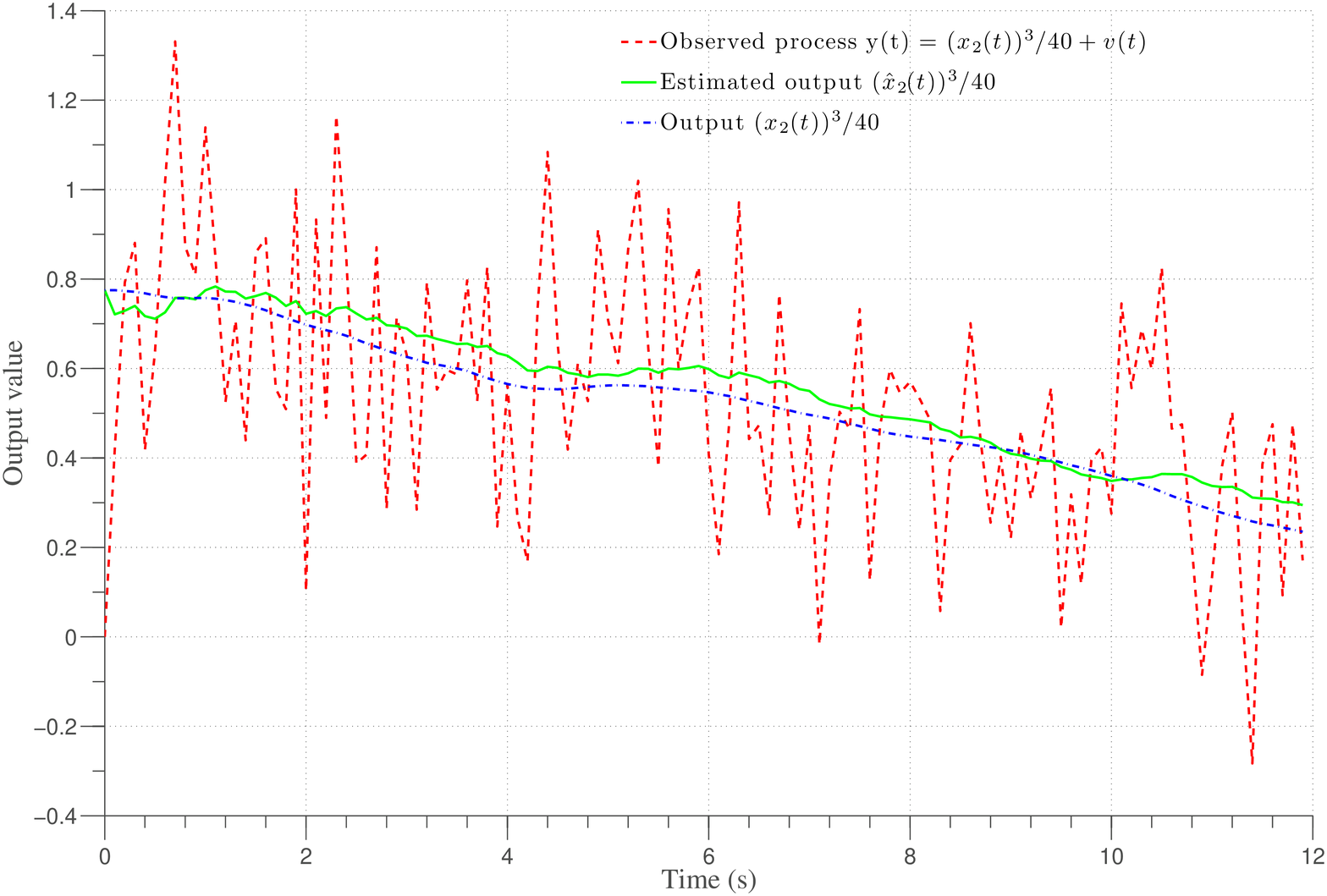}
\caption{Filtered measurement}
\label{\figname examplesimmeas}
\end{center}
\end{figure}

3. To reduce the dimensionality of the growth in basis terms
we use a k-means based clustering  approach  \cite{mackay2003information}
where the quadratics are clustered based on the locations of their argmin 
(ref. Sec.~\ref{\secname subsecpruning}). The simulation results of this filter design approach for the example considered are as shown in 
Fig.~\ref{\figname examplesimstatefilt}, 
\ref{\figname examplesimmeas}. Note that although each window is only $2$ units in width, the use of the repeated window generation
(a sliding window) leads to a filter capable of handling large measurement noises. 
 

\section{Conclusion and Future Directions} \label{\secname conclusions}
The details of the min-plus approach  described herein generate  filters for systems with nonlinear dynamics and nonlinear output. Its main 
novelty is in the  utilization  of the min-plus basis expansion of the value function coupled with the exploitation of
the linearity of the dynamic programming operator over such a (semi)-field. The salient features of this approach as discussed herein,
provide a structure which maybe used for a variety of applications. Further, the example serves as a template for future algorithmic developments.
A few  of the avenues along which a study of the  different stages used in these methods may be pursued are: (1) generating optimal min-plus fitting 
techniques. The approach to  fitting described herein provides one method (albeit not the optimal one); (2) the development of 
more sophisticated projection techniques for managing the growth in dimensionality (and an error analysis thereof); (3) the study of methods to 
speed up these stages. For instance some of the basis expansions may be done offline thereby helping speed up the real time operation of
the implementation. The application of these methods to the real time estimation of signals in various domains is a potentially fruitful 
theme for future research. 

\section*{Acknowledgment}   
The author would like to thank William \mbox{McEneaney} for helpful discussions and insights.

\bibliographystyle{unsrt}
\balance

%


\end{document}